\documentclass{amsart}
\usepackage{amssymb, amsmath}

\def\largerightarrow{   -\negthinspace\negthinspace -\negthinspace\negthinspace
                                       \negthinspace\longrightarrow }

\begin{document}
\title{}
\vbox{\hfil {\Large\bf THE DIRICHLET PROBLEM AND SPECTRAL  }\hfil }
\vbox{\hfil {\Large\bf   THEORY OF OPERATOR ALGEBRAS}\hfil}
\author{U. Haag}
\date{\today}
\maketitle
\noindent
{\bf 0. Introduction.}
\par\bigskip\noindent
Let $\, x\, $ be an element of a Banach algebra $\, A\, $ and $\,  U \supset sp_A ( x )\, $ any open neighbourhood of its spectrum. If $\, \mathfrak H ( U )\, $ denotes the Banach algebra of functions bounded and holomorphic in $\, U\, $ equipped with the supremum norm there is defined a holomorphic function calculus 
$$ {\omega }_{U , x} :\> \mathfrak H ( \mathcal U )\> \longrightarrow\> A $$
which is a unital homomorphism of Banach algebras and maps the identical function $\, id_z\, $ to 
$\, x\, $, cf. \cite{K-R}. Assuming that the holomorphic function calculus is contractive for each such neighbourhood $\, U\supset sp_A ( x )\, $, for example if $\, x\, $ is normal in some enveloping $C^*$-algebra, one then has a well defined contractive holomorphic function calculus 
$$ {\omega }_{A , x} :\> \mathfrak H ( sp_A ( x ) )\>\longrightarrow\> A $$
where $\, \mathfrak H ( sp_A ( x ) )\, $ denotes the uniform closure in $\, C ( sp_A ( x ) )\, $ of restrictions of functions holomorphic in some neighbourhood of $\, sp_A ( x )\, $.   
The purpose of this treatise is to extend the holomorphic function calculus in the case where $\, A\, $ is a unital operator algebra, which implies that it has a uniquely determined enveloping operator system, and assuming that the holomorphic function calculus is contractive, to a so called {\it harmonic function calculus} (recall that the real part of a holomorphic function is a harmonic function with uniquely determined conjugate harmonic function, up to a constant, cf. \cite{Be-So}) by considering the unique selfadjoint extension of the holomorphic function calculus to the enveloping operator systems. 
The point here is that although the extension of the holomorphic function calculus is only linear (and completely positive) its source can be shown to be linearly isometric with a commutative $C^*$-algebra by the main theorem below which to some extent allows a "continuous" function calculus even if not in the form of a $*$-homomorphism. 
In the course of this construction we prove a generalization of the well known 
Dirichlet problem to arbitrary compact subsets of the complex plane replacing the topological boundary by the Shilov boundary of the algebra of functions holomorphic in some neighbourhood of the compact set $\, X\, $ viewed as a subalgebra of $\, C ( X )\, $. The interplay between harmonic functions in the complex plane and holomorphic function theory enters with great subtlety in certain places and classical counterexamples, e.g. {\it l'\'{e}pine de Lesbesgue}, show that corresponding results do not hold in full generality in higher dimensions. 
We also consider real operator spaces and their injective envelopes since these arise quite naturally as subspaces of hermitian holomorphic functions on a domain which is symmetric with respect to reflection at the real axis. The results are then applied to define superpositive absolute values for certain hermitian elements of a super $C^*$-algebra $\, A\, $ in its enveloping graded operator system (see \cite{Ha1} for the definitions) and (superpositive) square roots in case of superpositive elements. Similarly one may define positive absolute values of certain elements of an operator algebra in its enveloping operator system. Finally we propose an analogue for the algebra of holomorphic functions replacing the Cauchy-Riemann equations by a set of partial differential equations of similar kind in higher dimensions which we believe is new, at least this was priorily unknown to the author.
\par\bigskip\bigskip\noindent
{\bf 1. The $*$-holomorphic functional calculus and real operator spaces.}
\par\bigskip\noindent
Assume given a hermitian element 
$\, x = x^*\, $ in a super $C^*$-algebra $\, A\, $ such that the closed commutative  subalgebra $\, A_x\, $ generated by $\, 1\, $ and $\, x\, $ embeds into a commutative $C^*$-algebra $\, C_x\, $, i.e. if $\, x\, $ is normal in the enveloping minimal $C^*$-algebra of $\, A_x\, $. Then the holomorphic function calculus 
$\, {\mathfrak H}_e ( sp ( x ) ) \rightarrow A_x\, $, where $\, {\mathfrak H}_e ( sp ( x ) )\, $ denotes the closure of the space of restrictions of entire holomorphic functions to the spectrum of $\, x\, $, sending the identity map $\, id_z : z \mapsto z\, $ to $\, x\, $ and $\, e ( z )\equiv 1\, $ to $\, 1\, $ is completely isometric and onto. Here 
$\, sp ( x )\, $ denotes the spectrum of 
$\, x\, $ in $\, A_x\, $ but one may also take the spectrum in any enveloping commutative  $C^*$-algebra for example $\, C_x\, $.
Since $\, x\, $ is hermitian its spectrum is symmetric with respect to reflection at the real line. Then 
$\, {\mathfrak H}_e ( sp ( x ) )\, $ is naturally a super $C^*$-algebra with involution determined by 
antilinear and multiplicative continuation of $\, id_z^* = id_z\, $ as follows from the "Small Reflection Theorem" of Schwarz (cf. \cite{Be-So}). It is immediate that the holomorphic function calculus is 
$*$-linear for this involution. Being a complete isometry the function calculus extends uniquely to a (bijective) contractive, hence 
graded isometric map of the enveloping graded operator systems $\, {\widehat{\mathfrak H}}_e ( sp ( x ) ) \rightarrow {\widehat A}_x\, $. If $\, x\, $ is superpositive, then $\, \iota ( x )\, $ is positive  in 
$\, {\widehat A}_x\, $ (see \cite{Ha1} for the definition of $\,\iota\, $ and other related notation), which implies that any point $\, z\, $ in its spectrum satisfies 
$\, Re ( z ) \geq 0\, $ and $\, \vert Im ( z )\vert \leq Re ( z )\, $, because the graded decomposition
$\, x = x_0 + x_1\, $ coincides with the decomposition into the selfadjoint part  $\, x_0\, $ and an 
antiselfadjoint part $\, x_1\, $. By the completely contractive (and hence completely superpositive) Gelfand transform 
$$ {\mathfrak G}_x :\> A_x\> \longrightarrow\> C ( sp ( x ) ) $$
this remains true even if $\, A_x\, $ does not embed into a commutative $C^*$-algebra and the holomorphic function calculus $\, {\mathfrak H}_e ( U ) \rightarrow A_x\, $ which is sort of an inverse to 
$\, {\mathfrak G}_x\, $ is not contractive for a given small open neighbourhood 
$\, U\supset sp ( x )\, $. 
Choose $\, \epsilon > 0\, $ arbitrarily small. By a wellknown result 
(cf. \cite{F} , \cite{Bo-St}) the element $\, id_z + \epsilon 1\, $ has a square root in 
$\, \mathfrak H ( sp ( x ) )\, $ (which again can be chosen superpositive and then is unique) so that 
for $\, \epsilon \to 0\, $ these elements converge towards the unique superpositive square root of 
$\, id_z\, $. The image of $\, \sqrt{ id_z}\, $ under holomorphic function calculus $\, {\omega }_{A , x}\, $ if defined is a square root for $\, x\, $ which is superpositive if the holomorphic function calculus is contractive (hence superpositive). This is certainly the case if $\, x\, $ is a (locally) normal element in the diction below.
As an example consider the super $C^*$-algebra $\, A_0\, $ of 
\cite{Ha1} generated by the single hermitian element of square zero 
$$  \begin{pmatrix} \>\>\> 1 &\>\>\>  1 \\ -1 & -1 \end{pmatrix}\> \in\> {\widehat M}_2 ( \mathbb C ) $$
with off-diagonal grading on $\, {\widehat M}_2 ( \mathbb C )\, $. One checks that the square root of the superpositive element 
$$ x\> =\> \begin{pmatrix} \>\>\> 3 & 1 \\ -1 & 1 \end{pmatrix} $$
is given by 
$$ \sqrt x\> =\> { 1\over 2\sqrt 2} \begin{pmatrix} \>\>\> 5 & 1 \\ -1 & 3 \end{pmatrix} \> . $$
which again is superpositive. The problem here is that the holomorphic function calculus in an open neighbourhood of $\, sp ( x )\, $ in the superpositive segment need not be contractive since the radical of $\, A_x = A_0\, $ is nontrivial so that the Gelfand transform $\, {\mathfrak G}_x\, $ is not isometric. Then the holomorphic function calculus in the superpositive segment cannot be expected to be superpositive. There seems to be no good reason why a superpositive element should always have a superpositive square root, though a counterexample still has to be found, i.e. we cannot give a decisive answer to this question here.
\par\smallskip\noindent
To the contrary, even in the normal case, the square of a superpositive element need not be superpositive. This means that one cannot define a superpositive absolute value of a hermitian element $\, x\, $ by the formula 
$\, \vert x\vert = \sqrt{ x^2}\, $ as in the case of a selfadjoint element. Instead one has to employ the 
graded function calculus in an enveloping graded $C^*$-algebra of \cite{Ha1} to obtain a reasonable definition. We will see that it is possible to define a superpositive "absolute value" for any hermitian element $\, x \, $ of a super-$C^*$-algebra $\, A\, $ {\it in the enveloping graded operator system of 
$\, A\, $} using the harmonic function calculus which is the unique positive extension of the holomorphic function calculus in a given domain $\, U\, $  containing the spectrum of 
$\, x\, $ and chosen large enough so that the holomorphic function calculus is contractive. The result will depend on the chosen domain, but in certain cases there is an "optimal" choice corresponding to a minimal domain. We will briefly sketch the construction in the case where the domain is a closed disk 
$\, D\, $ with radius equal to $\, \Vert x\Vert\, $. Consider the super $C^*$-algebra 
$\, \mathfrak H ( D ) = {\mathfrak H}_e ( D )\, $ of continuous functions holomorphic in the interior of 
$\, D\, $ with enveloping (commutative) graded $C^*$-algebra $\, C ( D )\, $ of continuous functions on 
$\, D\, $. Applying graded function calculus with respect to the graded absolute value on replacing the ordinary product by the $*$-product of \cite{Ha1} one finds that the result, for the identity map 
$\, x \bumpeq id_z\, $ is given by the continuous function $\, \vert x {\vert }_s\, $ such that 
$$ \vert x {\vert }_{s , 0}\> =\> \vert x_0 \vert \vee \vert x_1 \vert\> ,\quad 
\vert x {\vert }_{s , 1}\> =\> \pm i\, ( \vert x_0 \vert \wedge \vert x_1 \vert ) $$
where the absolute values on the right are the ordinary absolute values of the corresponding functions and $\, \vee\, $ denotes the maximum of two given realvalued functions, while $\, \wedge\, $ denotes the minimum of two realvalued functions, and the choice of sign in $\, \vert x {\vert }_1\, $ depends on whether the point under consideration lies in the upper halfplane (has positive imaginary part) or in the lower halfplane (has negative imaginary part). This definition leads to an overall continuous function, since the imaginary part vanishes on the real line. Consider its image in the Shilov boundary of 
$\,\mathfrak H ( D )\, $ which can be identified with the continuous functions on the boundary 
$\, \delta D\, $ of the disk. It is well known (and in any case easy to see) that $\, C ( \delta D )\, $ is 
linearly generated by the image of $\, \mathfrak H ( D )\, $ and its adjoint space, i.e. 
$\, \widehat{\mathfrak H} ( D ) \simeq C ( \delta D )\, $ in this particular case. We want to show that the same holds with $\, D\, $ replaced by an arbitrary compact subset $\, X \subseteq \mathbb C\, $, replacing $\, \mathfrak H ( D )\, $ by the uniform closure $\, \mathfrak H ( X )\, $  in $\, C ( X )\, $ of functions holomorphic in some neighbourhood of $\, X\, $. To get this one increases the complexity of 
$\, X\, $ step by step. 
\par\smallskip\noindent
Before doing so we insert a short discussion of a notion of injective envelope for a real operator space, i.e. a real closed subspace $\, \mathfrak R\subseteq \mathcal B ( \mathcal H )\, $. Consider 
$\, \mathcal B ( \mathcal H ) \simeq \mathcal B ( {\mathcal H}_{\mathbb R} ) \otimes \mathbb C\, $ as the complexification of the real $C^*$-algebra of bounded operators on a real Hilbert space.
Any element $\, x\in \mathfrak R\, $ can be written in the form $\, x = x_0 + i x_1\, $ with 
$\, x_0\, ,\, x_1 \in \mathcal B ( {\mathcal H}_{\mathbb R} )\, $. 
Then there is the real operator space $\, \overline{\mathfrak R}\, $, the image of $\,\mathfrak R\, $ by complex conjugation, which is treated as a $ {\mathbb Z}_2$-grading and is checked to be completely isometric, but we do not assume $\,\mathfrak R = \overline{\mathfrak R}\, $, so that 
the homogenous parts $\, x_0\, $ and $\, i x_1\, $ as above need not be contained in 
$\, \mathfrak R\, $. Replacing 
$\, \mathfrak R\, $ by the set of (conjugate) diagonal elements $\, \{ x \oplus \overline x \}\, $ in 
$\, \mathfrak R \oplus \overline{\mathfrak R}\, $ one gets that $\,\mathfrak R\, $ is completely isometric to the subspace 
$$  \widehat{\mathfrak R}\, =\> \left\{\, \begin{pmatrix} \>\>\> x_0 &  x_1 \\  - x_1 &  x_0 \end{pmatrix}\,\Bigm\vert\, x\in\mathfrak R\, \right\}\, \subseteq  \mathcal B ( {\mathcal H}_{\mathbb R} \oplus {\mathcal H}_{\mathbb R} )  $$
which can be seen by first rotating $\, x\oplus \overline x\, $ to 
$$ \begin{pmatrix} x_0 & i x_1 \\ i x_1 & x_0 \end{pmatrix} $$
and then conjugating with the unitary 
$$ \kappa\> =\> \begin{pmatrix} 1 &  0 \\ 0 & i \end{pmatrix} \qquad \> . $$
If $\, \mathfrak R\, $ happens to be a real operator system, i.e. a real subspace of $\, \mathcal B ( \mathcal H )\, $ which is selfadjoint and contains the identity element, then the same will be true for 
$\, \widehat{\mathfrak R}\, $. Namely, the complete isometry of $\, \mathfrak R\, $ and 
$\, \widehat{\mathfrak R}\, $ as above extends to a $*$-linear completely isometric real linear map
$$ \mathcal B ( \mathcal H )\, \longrightarrow\, M_2 ( \mathcal B ( {\mathcal H}_{\mathbb R} ) ) \> , $$
and by injectivity of $\, \mathcal B ( {\mathcal H}_{\mathbb R} )\, $ in the category of real operator spaces and real linear completely contractive maps, one sees that also the image of 
$\, \mathcal B ( \mathcal H )\, $ under this map is a real injective subspace since there exists an obvious completely contractive projection of $\, M_2 ( \mathcal B ( {\mathcal H}_{\mathbb R} ) )\, $ onto this image. Then one defines the injective envelope $\, I ( \mathfrak R )\, $  of the concrete real operator space $\, \mathfrak R\subseteq \mathcal B ( \mathcal H )\, $ to be the range of a minimal projection
$$ {\varphi }_{\mathfrak R} :  M_2 ( \mathcal B ( {\mathcal H}_{\mathbb R} ) ) \> \longrightarrow\> 
M_2 ( \mathcal B ( {\mathcal H}_{\mathbb R} ) ) \> , $$
i.e. $\, {\varphi }_{\mathfrak R}^2 = {\varphi }_{\mathfrak R}\, $, which is the identity restricted to 
$\, \widehat{\mathfrak R}\, $. Of course from what is said above the injective envelope 
$\, I ( \mathfrak R )\, $ may also be assumed to be a real subspace of 
$\, \mathcal B ( \mathcal H )\, $. One gets analogous results as in the complex case, including a {\it rigidity} property stating that any real linear completely contractive extension $\, \iota : I ( \mathfrak R ) \rightarrow I ( \mathfrak R )\, $ of the identity map of 
$\, {\widehat{\mathfrak R}}\, $ is the identity map of $\, I ( \mathfrak R )\, $ and by injectivity of 
$\, \mathcal B ( {\mathcal H}_{\mathbb R} )\, $ in the category of real operator spaces and completely contractive maps there is a corresponding {\it extension} (or injectivity) property for
$\, I ( \mathfrak R )\, $. Assume that $\, \mathfrak R\, $ is a real operator system. Then the map 
$\, {\varphi }_{\mathfrak R}\, $ can be taken to be selfadjoint, by a simple averaging process. Since it is also unital and completely contractive it is necessarily completely positive (compare the proof of Corollary 5.1.2 of \cite{E-R}).  One may define a product on $\, I ( \mathfrak R )\, $ by the formula 
$$ x * y\> =\> {\varphi }_{\mathfrak R} ( x\, y ) $$
if $\, x\, ,\, y\in I ( \mathfrak R )\, $.  
To get associativity of the product one uses that complete positivity of $\, {\varphi }_{\mathfrak R}\, $ entails, via  Stinesprings dilation theorem, a Schwarz inequality for the map and all of its amplifications. Then a standard argument as in the complex case gives associativity of the product as above. In addition one derives the $C^*$-condition 
$$ \Vert x * x^* \Vert = \Vert x {\Vert }^2 \> , \quad x\in I ( \mathfrak R ) $$
analogous to the complex case, so that $\, I ( \mathfrak R )\, $ admits via the GNS-construction a  
$*$-linear representation as a concrete real $C^*$-algebra on some real Hilbert space. There is a faithful and completely isometric representation, since 
for every selfadjoint element $\, x\, $ of the original real operator system $\, I ( \mathfrak R )\, $ there is a state $\,\rho\, $ with $\, \vert \rho ( x ) \vert = \Vert x \Vert\, $. Then its complexification is a complex $C^*$-algebra. 
\par\bigskip\bigskip\bigskip\bigskip\bigskip\bigskip\noindent
{\bf 2. The Dirichlet problem and the harmonic function calculus.}
\par\bigskip\noindent
We now return to our discussion of the Shilov boundary $\, \delta X\, $ of 
$\, \mathfrak H ( X )\, $ with $\, X\subset \mathbb C\, $ a compact subset in order to show that 
$\, \widehat{\mathfrak H} ( X ) \simeq C ( \delta X )\, $. There are three different cases of interest corresponding to different subspaces of holomorphic functions: the smallest is given by the closure of the space of restrictions of entire holomorphic functions to $\, X\, $ which is denoted 
$\, {\mathfrak H}_e ( X )\, $, the second is the closure of the space of functions holomorphic in some open 
neighbourhood of $\, X\, $ which is $\, \mathfrak H ( X )\, $.  For each of these spaces (unital operator algebras) there is a uniquely determined (in the sense of operator space theory) enveloping operator system consisting of a certain subspace of harmonic functions. These are denoted $\, {\widehat{\mathfrak H}}_e ( X )\, $ and $\, \widehat{\mathfrak H} ( X )\, $ respectively. Also for each of these spaces their Shilov boundary $\, {\delta }_e X\, $ resp. $\, \delta X\, $ can be realized as a compact subset of 
$\, X\, $. If $\, X\, $ is the closure of an open domain these are contained in the ordinary boundary of 
$\, X\, $ which for reasons of distinction is denoted 
$\, \dot\delta X = X\backslash \Breve X\, $ as follows from the maximum principle. 
The goal is to show that $\, {\widehat{\mathfrak H}}_e ( X ) \simeq C ( {\delta }_e X )\, $ and
$\, \widehat{\mathfrak H} ( X )\simeq C ( \delta X )\, $ for any compact subset $\, X\, $ of the complex plane. 
An immediate observation is that
if the result $\, \widehat{\mathfrak H} ( X_j ) \simeq C ( \delta X_j )\, $ is true for any two disjoint compact sets $\, X_1\, $ and $\, X_2\, $, then also for their union 
$\, X = X _1 \cup X_2\, $. This is easily checked from the fact that 
$\, \delta X = \delta X_1 \cup \delta X_2\, $ whence 
$\, C ( \delta X ) \simeq C ( \delta X_1 ) \oplus C ( \delta X_2 )\simeq 
\widehat{\mathfrak H} ( X_1 \cup X_2 )\, $. 
In the entire case things are seemingly much more complicated. 
So for $\, \delta X\, $ at least  the question is reduced in some sense to approximating an arbitrary compact set by sets which are finite disjoint unions of connected compact sets of a simpler type. 
For $\, {\delta }_e X\, $ the situation turns out to be quite similar. 
In principle the question is closely related to the wellknown 
Dirichlet Problem of Potential Theory, restricted to the complex plane, but with the additional requirement that the space of harmonic functions on a domain in question should be equal to the enveloping operator system of $\, \mathfrak H ( X )\, $. In the following the word {\it domain} will be used to mean a {\it bounded and connected} open subset of the complex plane which is equal to the interior of its closure and such that its boundary consists of a finite number of closed curves which are smooth except in finitely many places.
If $\, \Breve X \, $ is such a domain with compact closure $\, X\, $ then it is clear from the maximum principle that the Shilov boundary $\,\delta X\, $ of 
$\, \mathfrak H ( X )\, $ can be realized as a compact subset of the boundary 
$\, X \backslash \Breve X\, $. That $\,\delta X\, $ in fact consists of the whole boundary can be seen in the following way. At each smooth boundary point $\, z_0\in \dot\delta X\, $ there exists a small circle 
$\, \mathfrak C\, $ containing $\, z_0\, $ and such that the interior of the circumscribed disk has empty intersection with $\, X\, $. By a simple translation which is holomorphic one can assume that the center of the disk is given by the origin of the complex plane. Then there exists another circle $\,\mathfrak D\, $ around the origin such that $\, X\, $ is in the annulus with boundary $\, \mathfrak C\cup \mathfrak D\, $. The holomorphic transformation 
$\, \mathbb C\backslash \{ 0 \} \rightarrow \mathbb C \backslash \{ 0 \} \, ,\, z\mapsto z^{-1}\, $ sends the circle $\,\mathfrak C\, $ to the circle $\, {\mathfrak C}^{-1}\, $ such that the point $\, z_0^{-1}\, $ 
becomes a "convex" point of the transformed domain $\, X^{-1}\, $, meaning that there exists a closed disk (circumscribed by $\, {\mathfrak C}^{-1}\, $) containing the domain $\, X^{-1}\, $ such that its boundary contains the point corresponding to $\, z_0^{-1}\, $. Assuming as we may that this point lies on the positive real axis one sees that there exists a holomorphic function, for example $\, f ( z ) = e^z\, $ will do, which takes its maximum precisely in this point. Therefore $\, z_0\in\delta X\, $. Since $\, \delta X\, $ is closed it must agree with the usual boundary $\,\dot\delta X\, $. 
Let us denote the space of holomorphic (harmonic) functions in $\,\Breve X\, $ which admit a continuous extension to the boundary by 
$\, {\mathfrak H}_c ( X )\, $ (resp. $\, {\widehat{\mathfrak H}}_c ( X )\, $). 
The solution of the Dirichlet Problem asserts that $\, {\widehat{\mathfrak H}}_c ( X ) \simeq 
C ( \dot\delta X )\, $. Assuming that $\, \delta X = \dot\delta X =  X \backslash\Breve X\, $ which as seen applies to the type of domains we consider here one still has to show 
$\, \widehat{\mathfrak H} ( X ) = {\widehat{\mathfrak H}}_c ( X )\, $.
Let us call the former assertion I (ordinary Dirichlet problem), and the latter assertion II. These assertions are both contained in the statement $\, \widehat{\mathfrak H} ( X ) \simeq C ( \delta X )\, $ which makes sense for any compact set $\, X\subseteq\mathbb C\, $. 
In addition we consider assertion III, which is 
$\, {\widehat{\mathfrak H}}_e ( X ) = C ( {\delta }_e X )\, $. 
To keep in touch with holomorphic functions one may also ask: is 
$\, \mathfrak H ( X )\, $ the same as $\, {\mathfrak H}_c ( X )\, $ (assertion IV), and is 
(assertion V) $\, {\widehat{\mathfrak H}}_c ( X )\, $ the enveloping operator system of 
$\, {\mathfrak H}_c ( X )\, $ ? A positive answer to assertion II will imply assertion V (but not necessarily assertion IV), conversely knowledge of assertions IV and V implies assertion II. That these questions are not quite trivial can be seen from simple examples in the unit disk where there exist functions which are harmonic in the interior and extend continuously to the boundary, but such that their conjugate functions do not have a continuous extension (cf. \cite{Be-So}, p. 175). In a nice convex domain $\,\Delta\, $ without any holes one can find for every 
$\, \epsilon > 0\, $ a $\,\delta > 0\, $ such that rescaling the complex plane by the factor $\, 1 - \delta\, $ (assuming that the zero point is contained in $\, \Delta\, $) the corresponding harmonic function 
$\, h_{\delta}\, $ by rescaling of $\, h\, $ is $\epsilon $-close to 
$\, h\, $ on $\, \Delta\, $ and such that its conjugate is bounded and harmonic, hence continuous on 
$\,\delta\Delta\, $, but in a more general domain this is not seen very easily (as far as the author is concerned). For a convex domain one also checks that 
$\, {\delta }_e \Delta = \delta\Delta\, $ since any "convex" point $\, z_0\, $ of the boundary $\,\delta X\, $ is certainly in $\, {\delta }_e X\, $ and the boundary of a strictly convex domain consists solely of convex points. The case of a general convex domain follows from the "cut-off" procedure as described below.
We think it not unworthy to reprove known results on 
solutions of the Dirichlet problem fitting into our scheme using methods of Operator Space Theory. Thus we pretend to know nothing about such solutions beforehand except in the case of a closed disk where as mentioned above the statement is fairly obvious. The idea of proof is to establish several 
"moves" by which the statement for a given space is deduced from its validity for another related space or a certain class of spaces. In fact, instead of merely exhibiting stability of assertion I under each of these moves, we will also need to incorporate assertions II and III into the argument. It is convenient in certain arguments to assume these conditions from the start for a given domain in question. This property (I--III) is then shown to be stable by certain of the moves given below. As a starting point one verifies the stronger property for the unit disk $\, D\, $ of radius $1$ (centered without loss of generality in the origin). We have already seen that 
$\, {\widehat{\mathfrak H}}_c ( \Delta ) = \widehat{\mathfrak H} ( \Delta )\, $ holds for any convex domain. Then we only need to prove the other identity. If $\, f\, $ is any holomorphic function defined in some neighbourhood of $\, D\, $ then there exists a disk $\, D_R\, $ of radius 
$\, R > 1\, $ having the same center as $\, D\, $ where $\, f\, $ can be represented as a power series 
$$ f ( z )\> =\> \sum_{n = 0}^\infty\, a_n z^n  $$ 
such that the series of coefficients 
converges absolutely $\, \sum \vert a_n \vert < \infty\, $. This implies that $\, f\, $ can be uniformly approximated by holomorphic polynomials in $\, D\, $, i.e. by entire holomorphic functions.  
From the rescaling argument given above the result then follows for a disk of arbitrary radius by the fact that rescaling an entire harmonic function yields an entire harmonic function. 
We proceed by describing the first of the moves mentioned above. 
Assume that $\, X\, $ is a compact subset of the complex plane which is the closure of a connected domain $\, \Breve X\, $. Suppose that 
$\, X\, $ is divided by a straight line which without loss of generality may be assumed to be the real line on applying a translation (by a complex number) and a rotation (multiplication by a complex number of modulus one) both of which are holomorphic transformations, shifting $\, X\, $ into an appropriate position. Let $\, X^+\, $ be the intersection with the upper halfplane and $\, X^-\, $ its intersection with the lower halfplane. Assume that the reflection 
$\, \alpha \, $ at the real line takes $\,  X^-\, $ into 
$\, X^+\, $.  The first step is to show that if $\, X\, $ satisfies the strengthened version of assertions I--III, i.e. in particular $\, {\delta }_e X = \delta X\, $, and if
$\, X^0\, $ is the space obtained from $\, X\, $ by cutting off any finite number of components of $\, X^-\, $ (a component is supposed to mean a maximal connected subset of 
$\, X^- \cap \{ \, z\,\vert\, Im\, z < 0\, \}\, $, it then suffices to assume that the reflection map takes any of these components into $\, X^+\, $), then $\, X^0\, $ has the same property. 
Building on this result the next step shows that in case $\, X^0 = X^+\, $ the assertions I--III also hold for the symmetrization $\, \overline X = X + \alpha ( X )\, $ with respect to reflection at the real line. 
$\, \mathfrak H ( \overline X )\, $ is a super-$C^*$-algebra with involution 
$\, f^* ( z ) = \overline{ f ( \alpha ( z ))}\, $. Considering the real subspace of hermitian elements 
$\, \mathfrak H ( \overline X )^h\, $, which only take real values on the real line, it is easily seen to be completely isometric with its image space in $\, \mathfrak H ( X^+ )\, $ by restriction. Conversely, any element of $\, {\mathfrak H}_c ( X^+ )\, $ which takes only real values on $\, \mathbb R\cap X\, $ has a unique isometric extension to a hermitian element of $\, {\mathfrak H}_c ( \overline X )\, $ as follows from the Small  Reflection Theorem of Schwarz. The corresponding real subspaces of $\, \mathfrak H ( X )\, $ and $\, \mathfrak H ( X^+ )\, $ under restriction are denoted $\, \mathfrak H ( X )^h\, $ and
$\,\mathfrak H ( X^+ )^h\, $ respectively. One finds that the natural graded extension of the maps 
$\,\mathfrak H ( \overline X )^h \leftrightarrow \mathfrak H ( X^+ )^h\, $ of real operator systems is still a complete isometry, the restriction map being completely contractive in any case and its inverse being induced by the flip of the complex plane and complex conjugation. We denote by 
$\, \widehat{\mathfrak H} ( X )^s\, $ and $\, \widehat{\mathfrak H} ( X^+ )^s\, $ the images of the completely isometric restriction of the even (symmetric) part of $\, \widehat{\mathfrak H} ( \overline X )\, $, and correspondingly by $\, \widehat{\mathfrak H} ( X )^{a}\, $ and $\, \widehat{\mathfrak H} ( X^+ )^{a}\, $ the images of the completely isometric restriction of the antisymmetric (odd) part of 
$\, \widehat{\mathfrak H} ( \overline X )\, $ to $\, X\, $ and $\, X^+\, $ respectively.
We also consider the space of holomorphic (or harmonic) functions on 
$\, {\Breve X}^c = \Breve X\backslash ( \overline{ X \backslash X^0} )\, $ with 
$\, \overline {X\backslash X^0}\, $ the symmetrization of $\, \Breve X\backslash {\Breve X}^0 \, $. Note that the restriction map 
to $\, \mathfrak H ( X^c )\, $ is still completely isometric for the subspace 
$\, \mathfrak H ( \overline X )^h\, $ since any element in this subspace takes its maximum on the boundary component of $\, \overline X\, $ which intersects the boundary of 
$\, X^+\, $ which is the same as $\, \delta \overline X \cap \delta X^c\, $.
Consider the super-$C^*$-algebra $\, {\mathfrak H}_e ( \overline X )\, $. Restricted to the real line the involution is simply complex conjugation, and the Stone-Weierstrass-Theorem implies that the image of 
$\, {\mathfrak H}_e ( \overline X )\, $ by evaluation on the real line is dense in 
$\, C ( \mathbb R \cap X )\, $ because it contains the constants, separates points, and is closed under complex conjugation. Then an arbitrary element of $\, f\in \widehat{\mathfrak H} ( X^0 )\, $ can be assumed to be uniformly small restricted to the real line, say to the order of some arbitrarily chosen 
$\, \epsilon > 0\, $ modulo the image of $\, {\mathfrak H}_e ( \overline X )\, $ which coincides with the complexification of $\, {\mathfrak H}_e ( X )^h\, $ (and $\, f \equiv 0\, $ on the intersection of the boundary of $\, X\, $ with the real line which only consists of finitely many points). 
Let $\, {\widehat{\mathfrak H}}_0 ( X )\, $ denote the subspace of 
$\, {\widehat{\mathfrak H}}_e ( X )\, $ whose elements are zero evaluated on 
$\, \delta X \cap \mathbb R\, $. Let $\,\epsilon > 0\, $ be given and choose a natural number 
$\, N\, $ such that $\, 1 / N < \epsilon / 2\, $. Given any continous function $\, f\, $ on 
$\, \delta X^0\cap \delta X\, $ which vanishes on the intersection with the real line and of norm one there exists a harmonic function $\, h\in {\widehat{\mathfrak H}}_0 ( X )\, $ with $\, \Vert h \Vert = 1\, $ such that $\, h\, $ agrees with $\, f\, $ on the boundary component $\, \delta X^0\cap\delta X\, $, by assumption on 
$\, X\, $. Consider the restriction of $\, h\, $ to $\, h^c\in \widehat{\mathfrak H} ( X^c )\, $.
Then $\, h^c\, $ is uniquely determined by the corresponding continuous function on $\,\delta X^c\, $ which is also denoted $\, h^c\, $, and identifies with $\, f\, $ on the intersection 
$\, \delta XÊ\cap \delta X^c\, $. There is a natural identification of $\, C_0 ( \delta X^c )\, $ with 
$\, C_0 ( \delta X )\, $ (the subspaces of continuous functions vanishing on 
$\, \delta X\cap\mathbb R\, $) given by flipping the boundary component of $\, X^c\, $ which does not intersect with $\,\delta X\, $ by reflecting it at the real line onto the corresponding boundary component of 
$\, X\, $, and taking the negative of the values of a function on this boundary component. Let $\, h_1\, $ denote the image of $\, h^c\, $ by this procedure. 
Then from the assumption on $\, X\, $ the function 
$\, h_1\, $ again defines an element of $\, {\widehat{\mathfrak H}}_0 ( X )\, $, and one may repeat this process by taking its restriction $\, h_1^c\, $ and so on, leading to a sequence of 
elements $\, \{ h_n {\}}_n\, $ of elements in $\, {\widehat{\mathfrak H}}_0 ( X )\, $. One notes that the values $\, h_n ( z ) \, $ remain constant for $\, z\in\delta X\cap \dot\delta X^c\, $ and that the norms are monotonously decreasing by this process
$$ \Vert h \Vert \geq \Vert h_1 \Vert \geq \cdots \geq \Vert h_n \Vert \geq \cdots \> . $$
Define 
$$ h^N\> =\> {1\over N}\, \sum_{k=1}^N\, h_k \> . $$
Then $\, \Vert  h^N_1 - h^N \Vert \leq 2 / N \Vert h \Vert < \epsilon\, $ which implies that the restriction of 
$\, h^N\, $ to $\, \widehat{\mathfrak H} ( \overline{ X\backslash X^0} )\, $, which space carries a natural superinvolution, is $\epsilon $-close to being antisymmetric, i.e. 
$$ \Vert h^N {\vert }_{\overline{X\backslash X^0}} + h^NÊ{\vert }_{\overline{X\backslash X^0}}\circ \alpha \Vert < \epsilon $$
In particular it must be $\epsilon $-close to zero on the real line. 
This shows that given any function $\, f\, $ and $\, \epsilon > 0\, $ as above there exists a harmonic function $\, h_{\epsilon }\in {\widehat{\mathfrak H}}_e ( X^0 )\, $ which agrees with $\, f\, $ on 
$\, \delta X^0\cap \delta X\, $ and such that 
$\, \Vert h_{\epsilon } {\vert }_{\delta X^0\cap\mathbb R} \Vert < \epsilon\, $. 
But then again by the argument above there exists an element  
$\, h_0\in {\widehat{\mathfrak H}}_0 ( \overline X )\, $ which is $\epsilon $-close to $\, f\, $ restricted to the real line, and by the previous argument again there exists another element 
$\, h_{0 , \epsilon }\in {\widehat{\mathfrak H}}_0 ( X^0 )\, $ which identifies with $\, h_0\, $ on 
$\, \delta X^0\cap \delta X\, $ and is $\epsilon $-close to zero on the real line so that the linear combination $\, h_f = h_{\epsilon } + h_0 - h_{0 , \epsilon }\, $ is $3\epsilon $-close to 
$\, f\, $ on $\,\delta X^0\, $. One concludes that $\, {\widehat{\mathfrak H}}_e ( X^0 )\, $ is dense in 
$\, C ( \delta X^0 )\, $. Since it is also closed assertions I--III follow. 
This argument suffices to show that assertions II and III are valid also for every convex domain. Indeed such a domain can be manufactured from a closed disk by an infinite number of "cut-offs" as above. In each finite step the corresponding space $\, {\Delta }_n\, $ will satisfy the  assertions I--III. Then it follows from the inductive limit argument given below that the limit space $\, \Delta = \cap_n\, {\Delta }_n\, $ also inherits the property 
$\, {\widehat{\mathfrak H}}_e ( \Delta ) = C ( {\delta }_e \Delta )\, $ and by convexity one has 
$\, {\delta }_e \Delta = \delta\Delta\, $ (for any point $\, z_0\, $ in the boundary of $\,\Delta\, $ there exists a triangle containing $\,\Delta\, $ such that its boundary contains $\, z_0\, $) which gives 
assertions I--III for $\,\Delta\, $.
The assumption that the components of $\, X^-\, $ to be cut off should be reflected into $\, X^+\, $ can be removed by "nibbling away" each component in several consecutive steps, such that in each single step the corresponding assumption holds for the smaller portions to be cut off. 
This result can then be used to derive a corresponding result also for the symmetrization
$\,\overline X\, $ of $\, X^+\, $ with respect to reflection at the real line. In fact since by the process above starting with an arbitrary element 
$\, h\in {\widehat{\mathfrak H}}_0 ( X )\, $ the sequence of functions $\, \{ h^NÊ\}\, $ converges uniformly to zero on the real line the limit function gives a well defined element 
$\, \overline h\,\in {\widehat{\mathfrak H}}_0 ( X^+ )\, $, so that the uniform convergence extends to all of $\, X^+\, $. Since the elements $\, h^N\, $ are increasingly close to being antisymmetric, the sequence converges uniformly in $\, X\, $, giving a well defined completely contractive projection 
$\, q : {\widehat{\mathfrak H}}_0 ( X ) \rightarrow {\widehat{\mathfrak H}}_0 ( X )\, $ whose range is precisely $\, \widehat{\mathfrak H} ( X )^{a}\, \simeq {\widehat{\mathfrak H}}_c ( \overline X )^{a}
= \widehat{\mathfrak H} ( \overline X )^{a}\, $ showing that every antisymmetric continuous function on 
$\, \delta \overline X\, $ occurs as the boundary value of some odd harmonic function and can be approximated by antisymmetric harmonic functions defined in some symmetric neighbourhood of 
$\,\overline X\, $. More effort is needed for the symmetric part of the statement. 
For technical reasons we assume that $\, \delta X\cap X^-\, $ is contained in a boundary segment of constant curvature of $\, \delta X\, $ such that the disk whose boundary contains this segment lies wholly in $\, X\, $. In particular $\, \mathbb R\cap\delta X\, $ consists of precisely two points which are both contained in the segment of constant curvature. This will be sufficient for our purposes. Choose any contracive selfadjoint linear map 
$\, \rho : C_0 ( \delta \overline X ) \rightarrow C_0 ( \delta X ) \simeq {\widehat{\mathfrak H}}_0 ( X )\, $ over the identity of $\, C_0 ( \delta X\cap\deltaÊ\overline X )\, $ (for example by trivial extension on the complementary segment).
Now consider the following process on 
$\, {\widehat{\mathfrak H}}_{0 , \mathbb R} ( X ) \simeq C_{0 , \mathbb R} ( \delta X )\, $ (the subspace of realvalued functions vanishing in $\, \mathbb R\cap \delta X\, $). 
Given an arbitrary element $\, h\, $ in this subspace consider the reflection of the boundary component 
$\, \delta X \cap X^-\, $ at the real line, and let $\, h^c\, $ denote the continuous function on 
$\, ( \delta X\cap X^+ ) \cup \alpha ( \delta X\cap X^- )\, $ given by evaluation of the harmonic function 
$\, h\, $ on this subset. Flipping the second component by reflecting it at the real line and transcribing the values for corresponding pairs of points defines another element 
$\, h_1\in C_{0 , \mathbb R} ( \delta X )\, $ whose values agree with those of $\, h\, $ on the boundary component $\, \delta X\cap X^+\, $, so the process may be iterated infinitely often leading to a sequence of elements $\, \{ h_k {\} }_k\subseteq \widehat{\mathfrak H} ( X )\, $ which are uniformly bounded by the norm of $\, h\, $. As above one defines 
$$ h^N\> =\> {1\over N}\, \sum_{k=1}^N\, h_k $$
and checks that each element $\, h^N\, $ is approximately symmetric to the order of $\, 2 / N\, $ on the symmetric subset $\, \overline{ X\cap X^+}\, $. Choosing a countable dense subset 
$\, \{ z_j \,\vert\, j\in\mathbb N \}\, $ of $\, X\, $ and a countable dense subsequence 
$\, \{ f_l\,\vert\, l\in \mathbb N \}\, $ in the unit ball of $\, C ( \delta X )\, $  there exists a subsequence 
$\, \{ N_k \} \subseteq\mathbb N\, $ such that for each pair of indices $\, ( j , l )\, $ the sequence 
$$  f^{N_k}_l ( z_j ) \> \buildrel {k\to\infty }\over\largerightarrow\> {\overline f}_l ( z_j ) $$ 
converges to a well defined limit value. 
The Poisson Integral Formula for harmonic functions in the interior of a disk then shows that the partial derivatives of an arbitrary harmonic function $\, h\, $ at a given point $\, z_0\, $ in the interior are bounded by a maximal value depending only on the norm of 
$\, h\, $ and the distance of the point $\, z_0\, $ to the boundary. Therefore one concludes that the sequence $\, ( f^{N_k}_l ( z ) {)}_k\, $ converges uniformly on compact subsets in the interior of $\, X\, $ to a limit function
$\, {\overline f}_l\, $ which is symmetric in (the interior of) $\, \overline {X\cap X^-}\, $ and therefore extends to a symmetric harmonic function in the interior of  $\, \overline X\, $. 
The Poisson Integral Formula also shows that for each point $\, z_0\, $ in the prescribed segment of constant curvature of the boundary not in $\,  \mathbb R\cap \delta X\, $ the limit function 
$\, {\overline f}_l ( z )\, $ converges to the boundary value $\, f_l ( z_0 )\, $ as $\, z\to z_0\, $, since the values of the functions 
$\, \{ f^{N_k}_l\} \, $ remain constant on $\  \delta X\cap X_+\, $ and uniformly bounded in the rest of the circle so that the directional derivatives  of $\, {\overline f}_l ( z )\, $ as $\, z\, $ tends to a point in the prescribed segment can be estimated against the directional derivatives of the original function 
$\, f_l ( z )\, $ (which without loss of generality are assumed to be bounded). Then it easily follows that the same holds for an arbitrary point $\, z_0\in ( \delta\overline X \backslash (\mathbb R\cap\delta X )\, $.
It remains to show that the values $\, \{ {\overline f}_l ( z ) \}\, $ converge to zero whenever $\, z\, $ converges to a boundary point of 
$\, \mathbb R\cap \delta X\, $. Let $\, {\overline g}_l\, $ be the unique (antisymmetric) conjugate harmonic function of $\, {\overline f}_l\, $ in the interior of $\, \overline X\, $ which vanishes on the intersection with 
$\,\mathbb R\, $, so that $\, {\overline h}_l = {\overline f}_l + i {\overline g}_l\, $ defines a hermitian holomorphic function in the interior of $\,\overline X\, $ with real part equal to $\, {\overline f}_l\, $. It then suffices to show that $\, {\overline h}_l ( z )\, $ converges to zero if 
$\, z \to z_0, z_0\in \mathbb R\cap\delta X\, $. Without loss of generality one may assume by the same reasoning as above by an appropriate choice of $\, f_l ( z )\, $ that the conjugate function 
$\, {\overline g}_l\, $ converges to an antisymmetric  continuous function on the boundary 
$\, \delta\overline X\, $ from Poisson's integral formula applied to a disk containing the boundary segment of constant curvature around $\, z_0\, $ as above.
One then checks that Cauchy's integral formula
$$ {\overline h}_l ( z )\> =\> {1\over 2\pi i} \int_{\delta \overline X} { {\overline h}_l ( \zeta )\over \zeta - z} d\zeta $$
applies to each point $\, z\, $ in the interior of $\,\overline X\, $ approximating this integral by a sequence of corresponding integrals over the boundary of some subset $\, \Breve X\subseteq \overline X\, $ by cutting out a small circle around each point of $\, \mathbb R\cap \delta X\, $. Then approaching $\, z_0\, $ from the open interval of $\, \mathbb R\, $ in the interior of $\,\overline X\, $ one may replace the boundary of $\,\overline X\, $ by the boundary $\, \mathfrak C = {\mathfrak C}_+ \cup {\mathfrak C}_-\, $ of the union of the two disks which are determined by the segments of constant curvature in this formula if $\, {\mathfrak C}_+\, $ denotes the boundary segment contained in $\, X_+\, $ and $\, {\mathfrak C}_-\, $ the segment in $\, X_-\, $.  Relating Cauchy's integral formula with Poisson's integral formula in either of the two disks one finds that the values of the Poisson integrals of both $\, {\overline f}_l \, $ as well as 
$\, {\overline g}_l\, $ over each of the curves $\, {\mathfrak C}_{\pm }\, $ must converge to zero 
for $\, t\to z_0\, ,\, t\in\mathbb R\cap\overline X\, $ since the corresponding boundary functions on 
$\, {\mathfrak C}_{\pm }\, $ can be extended continuously on the complementary segment of each disk by the zero function thus defining some harmonic function in the interior of the disk which converges to 
zero for $\, t\to z_0\, $ and such that the corresponding Poisson integral formula reduces to an integral over the nontrivial segment $\, {\mathfrak C}_{\pm }\, $. Cauchy's integral formula for $\, z \in \mathbb R\cap \overline X\, $ can be retrieved from 
the Poisson integral with respect to the boundary values of $\, {\overline f}_l\, $ and $\, {\overline g}_l\, $ along the curve $\, \mathfrak C\, $ showing that $\, {\overline h}_l ( t ) \to 0\, $ as $\, t\to z_0\, $. We leave the details of this construction to the reader. This then implies that $\, \{ f_l^{N_k} \}\, $ converges uniformly on $\, \mathbb R\cap\overline X\, $ to $\, {\overline f}_l\, $ with $\, {\overline f}_l{\vert }_{\mathbb R\cap\overline X} \, $ continuous. But then one must have uniform convergence in every point with continuous extension to the boundary. Since the sequence 
$\, \{ f_lÊ\}\, $ is dense in the unit ball of $\, C ( \delta X )\, $ this process defines a unital completely contractive and selfadjoint real linear map 
$$ \sigma  :  {\widehat{\mathfrak H}}_{0 , \mathbb R} ( X ) \> \longrightarrow\> 
{\widehat{\mathfrak H}}_c ( \overline X )^s \>  . $$
The composition 
$$ \sigma\circ\rho : \> C_0 ( \delta \overline X )^sÊ\longrightarrow \> C_0 ( \delta \overline X )^s \> $$
is just the identity map since the process defined above leaves invariant the boundary values on 
$\, \delta \overline X\cap X_+\, $. This shows that any symmetric continuous function on 
$\,\delta\overline X\, $ vanishing on the intersection with the real line admits an harmonic extension to the interior. It is then easy to separate the remaining two points in $\, \mathbb R\cap\delta\overline X\, $ from their complement by some nonvanishing harmonic functions.  
Putting together the hermitian and the antisymmetric part of the argument one concludes that 
$\, {\widehat{\mathfrak H}}_c ( \overline X )\, $ equals $\, C ( \delta \overline X )\, $ (assertion I).
Assertions II and III  will follow in general from the inductive limit procedure described below which asserts that if $\, X = \cap_n X_n\, $ and each $\, X_n\, $ satisfies 
$\, {\widehat{\mathfrak H}}_e ( X_n ) \simeq C ( {\delta }_e X_n )\, $ then the same holds for $\, X\, $ (or in case of assertion II the conclusion for $\, X\, $ also follows if only assertion I holds for all the $\, X_n\, $ since the limit of the spaces $\, {\widehat{\mathfrak H}}_c ( X_n )\, $ in $\, \widehat{\mathfrak H} ( X )\, $ by restriction is obviously the same as the limit of the spaces $\, \widehat{\mathfrak H} ( X_n )\, $). The 
$\, X_n\, $ can be taken to be finite unions of closed disks of the same diameter for which assertions I--III are proved below. 
\par\noindent
The next step consists in showing that if assertion I is true for a given subset 
$\, X_1\, $ which is the closure of a domain $\, \Breve X_1\, $ and if $\, X_2\, $ is a closed disk with nontrivial intersection $\, X_1 \cap X_2\, $ having finitely many (one or two is sufficient for our purposes) connected components with connected complement 
$\, X_2 \backslash X_1\, $, then assertion I also holds for $\, X = X_1 \cup X_2\, $. For technical reasons we again assume that in a small neighbourhood of each point of $\, \delta X_1 \cap \delta X_2\, $ the boundary $\, \delta X_1\, $ has constant curvature such that the prescribed disk is contained in 
$\, X_1\, $. Then from the previous steps above we may enlarge 
$\, X_2\, $ by adjoining a collar around the part of its boundary intersecting with $\, X_1\, $ (whose boundary is, say, a polygonial path approximately parallel to the original boundary and intersecting the boundary of $\, X_1\, $ in the prescribed segments as above. This enlarged "disk" is denoted $\, {\tilde X}_2\, $. The collar can be  made to fit the boundary segment of $\, \delta X_1\, $ which is inbetween the points constituting the intersection of $\, \delta X_1\, $ with the original disk boundary and the intersection with boundary of the enlarged space $\, {\tilde X}_2\, $ so that $\, \delta X_1\, $ and $\, \delta {\tilde X}_2\, $ both share two (or four if any) small segments of constant curvature, by manufacturing the collar from a (finite) number of symmetrizations and (possibly infinite number of) cut-offs along straight line segments (using the inductive limit procedure given below). One can avoid using the inductive limit procedure by assuming that the boundary segment of 
$\, X_1\, $ in question should have the same curvature as the boundary of the disk $\, X_2\, $ (for example if $\, X_1\, $ is a finite union of disks of constant diameter equal to the diameter of 
$\, X_2\, $). Then a simple symmetrization for each boundary segment in question will do and one gets
$\, {\widehat{\mathfrak H}}_c ( {\tilde X}_2 ) = C ( \delta {\tilde X}_2 )\, $. 
Consider the (pairs of) boundary points of the intersection of 
$\, \delta X_1\, $ with the boundary of the original disk and the corresponding boundary components of 
$\, \delta ( X_1\cap X_2 )\, $ joining the two points in each pair, the first being a boundary component of 
$\, \delta X_1\, $, and the second of $\,\delta X_2\, $ (in case that $\, \delta X_1 \cap \delta X_2\, $ is empty one gets two closed curves). Given a (realvalued) continuous function $\, f\in C ( \delta X )\, $ vanishing on $\, \delta X_1\cap \delta X_2\, $ the task is to find a function on $\, X\, $ which is harmonic in the interior and equals $\, f\, $ on the boundary. By assumption on $\, X_1\, $ the trivial extension of the restriction of $\, f\, $ to $\,\delta X_1\, $ defines an element $\, g\in {\widehat{\mathfrak H}}_0 ( X_1 )\, $ (the suffix $\, 0\, $ is used as above to denote functions vanishing in $\,\delta X_1\cap \delta X_2\, $), also the trivial extension of $\, f\, $ to 
$\, \delta X_2\, $ defines an element $\, h\in {\widehat{\mathfrak H}}_0 ( X_2 )\, $. One may evaluate the function $\, g\, $ on $\, \delta X_2\cap X_1\, $ and replace the trivial extension of $\, f\, $ on 
$\,\delta X_2\, $ by using the values of $\, g {\vert }_{\delta X_2\cap X_1}\, $ instead of the zero function. This defines another element $\, h_1\in {\widehat{\mathfrak H}}_0 ( X_2 )\, $. Similarly one may replace the trivial extension of $\, f\, $ on $\, \delta X_1\, $ by using the function 
$\, h {\vert }_{\delta X_1\cap X_2}\, $ instead. This defines an element 
$\, g_1\in {\widehat{\mathfrak H}}_0 ( X_1 )\, $. Repeating this process infinitely often leads to a sequence of pairs $\, ( g_n , h_n {) }_{n\in\mathbb N}\, $ in 
$\, {\widehat{\mathfrak H}}_0 ( X_1 ) \times {\widehat{\mathfrak H}}_0 ( X_2 )\, $. Since any harmonic function takes its maximum value on the boundary of the domain, the norms of this sequence are uniformly bounded by $\, \Vert f \Vert\, $. One may then consider the series 
$$ g^N\> =\> {1\over N} \sum_{k=1}^N\, g_k\> ,\quad h^N\> =\> {1\over N} \sum_{k=1}^N\, h_k \> , $$
which is also uniformly bounded and such that 
$$ \Vert g^N {\vert }_{X_1\cap X_2} - h^N {\vert }_{X_1\cap X_2} \Vert \leq 2 / N \> . $$
As above one shows the existence of a subsequence $\, ( g^{N_k} {\} }_k\, $ converging uniformly to a limit function in the interior of $\, X_1\, $ and near the boundary segment of constant curvature where the values of $\, \{ g^{N_k} \}\, $ remain constant, which together implies uniform convergence on the boundary segment $\, \delta {\tilde X}_2\cap X_1\, $. Extending the limit function by the constant values of $\, f\, $ on the remaining boundary segment of $\, \delta {\tilde X}_2\, $ shows uniform convergence of $\, \{ g^{N_k} \}\, $ on all of 
$\, \delta X_1\, $ which in turn implies uniform convergence of $\, \{ h^{N_k} \}\, $ on all of 
$\, \delta X_2\, $ such that the harmonic functions obtained in the limit agree on $\, X_1\cap X_2\, $, hence define an overall harmonic function $\, h_f\, $ with boundary values given by $\, f\, $. 
This proves assertion I. If $\, X\, $ is a finite union of closed disks of the same diameter then it can be approximated by a sequence of sets of the same form on choosing the diameter of the disks in the approximating spaces slightly larger. Then since each of these spaces satisfies assertion I as proved above the space $\, X\, $ will satisfy $\, \widehat{\mathfrak H} ( X ) \simeq C ( \delta X )\, $ from the inductive limit argument given below. One gets assertion III for such a space by the following 
scheme. Let $\, X_j\, ,\, j = 1 ,\cdots n\, $ be closed disks with $\, X = \cup_{j=1}^n\, X_j\, $ so that for each $\, j\, $ one has $\, {\mathfrak H}_e ( X_j ) = \mathfrak H ( X_j )\, $.
Let $\, x_j\in \mathcal B ( \mathcal H )\, $ be normal elements with $C^*$-spectrum equal to 
$\, X_j \subseteq sp ( x_j )\, $. Then from the fact that $\, {\widehat{\mathfrak H}}_e ( sp ( x_j ) ) \simeq {\widehat{\mathfrak H}}_e ( X_j )\simeq C ( \delta X_j )\, $ by the natural restriction map one easily sees that $\, sp ( x_j ) = X_j\, $ in each case.
Put $\, x = x_1 \oplus \cdots \oplus x_n\in\mathcal B ( \mathcal H\oplus\cdots \oplus \mathcal H )\, $. Then 
$\, x\, $ is normal with 
$\,  X \subseteq sp ( x )\, $ . For any open neighbourhood $\, U\supset sp ( x )\, $ the holomorphic function calculus 
$$ {\omega }_{U , x} :\> \mathfrak H ( U )\> \longrightarrow\> A_x $$
for functions holomorphic in $\, U\, $ is well defined and if $\, sp ( x ) \subset V\subseteq U\, $ is any smaller neighbourhood then the image of a holomorphic function $\, f_U\, $ in $\, A_x\, $ is the same as the image of its restriction $\, f_V\, $ under $\, {\omega }_{V , x}\, $. Let 
$\, \Breve{\mathfrak H} ( sp ( x ) )\subseteq \mathfrak H ( sp ( x ) )\, $ denote the dense subspace generated by  all restrictions of functions holomorphic in some neighbourhood of 
$\, sp ( x )\, $. This implies that one gets a well defined (though possibly unbounded) holomorphic function calculus
$$ {\omega }_x :\> \Breve{\mathfrak H} ( sp ( x ) )\> \longrightarrow\> A_x  $$
which is an isometry restricted to $\, {\mathfrak H}_e ( sp ( x ) )\, $ by normality of $\, x\, $. One also has a holomorphic function calculus 
$$ {\omega }_X :\> \mathfrak H ( X )\> \longrightarrow A_{x_1} \oplus\cdots \oplus A_{x_n} $$
factoring over the diagonal embedding $\, \mathfrak H ( X ) \subseteq \mathfrak H ( X_1 ) \oplus\cdots\oplus \mathfrak H ( X_n )\, $ by restriction which is contractive since each $\, X_j\, $ is a disk whence 
$\, \mathfrak H ( X_j ) = {\mathfrak H}_e ( X_j )\, $. Then it is also contractive on the restricted image of 
$\, \Breve{\mathfrak H} ( sp ( x ) )\, $ in $\, \mathfrak H ( X )\, $. By naturality the image of 
$\, \Breve{\mathfrak H} ( sp ( x ) )\, $ under $\, {\omega }_X\, $ is contained in $\, A_x\, $ viewed as a subalgebra of $\, A_{x_1} \oplus\cdots\oplus A_{x_n}\, $ via the diagonal embedding. This implies for one thing that 
$\, {\omega }_x\, $ must be contractive itself and hence extends to all of $\, \mathfrak H ( sp ( x ) )\, $. 
Composing $\, {\omega }_x\, $ with the (contractive) Gelfand transform 
$$ {\mathfrak G}_x :\> A_x\> \longrightarrow C ( sp ( x ) ) $$
with image equal to $\, {\mathfrak H}_e ( sp ( x ) )\, $ gives a contractive projection 
$$ \mathfrak H ( sp ( x ) )\> \longrightarrow {\mathfrak H}_e ( sp ( x ) )\> . $$
But from the very definition of the holomorphic function calculus evaluating the image of any 
holomorphic function $\, f_U\, $ under this composite map at an arbitrary point $\, z_0\in sp ( x )\, $  
gives back the value $\, f_U ( z_0 )\, $. Thus 
$\, \mathfrak H ( sp ( x ) ) = {\mathfrak H}_e ( sp ( x ) )\simeq {\mathfrak H}_e ( X )\, $. One may ask to what extent $\, sp ( x )\, $ is determined by $\, X\, $. The answer is that $\, sp ( x )\, $ is just the domain obtained by filling up all the "holes" in $\, X\, $. In particular every point $\, z_0\, $ of the "outer boundary" of $\, X\, $ for which there exists a continuous path $\, \{ z_t\,Ê\vert\, t\geq 0 \}\, $ connecting $\, z_0\, $ with the point at infinity and such that $\, z_t\notin X\, $ for $\, t > 0\, $, is also a point of the boundary of 
$\, sp ( x )\, $. Indeed if such a path is given then there exists a positive number $\, r > 0\, $ such that at each point $\, z_t\, $ one may attach a closed disk $\, D_t\, $ of radius $\, r\, $ with $\, z_t\in \delta D_t\, $ and no point in the interior of $\, D_t\, $ contained in $\, X\, $. From compactness of $\, sp ( x )\, $ there exists a minimal parameter $\, t_1\, $ such that the intersection of $\, sp ( x )\, $ with the tube provided by the union of the disks $\, \cup_t\, D_t\, $ is contained in $\, \cup_{t\leq t_1}\, D_t\, $. If $\, z_0\, $ is not in 
$\, \delta sp ( x )\, $ the intersection of $\, sp ( x )\, $ with the interior of $\, D_0\, $ must be nonempty so choosing $\, r\, $ small enough one can assume without loss of generality that $\, t_1 > 0\, $ and that 
$\, D_{t_1}\, $ has empty intersection with $\, X\, $. By minimality there must be a point $\, z_1\, $ of 
$\, sp ( x )\, $ contained in the boundary of $\, D_{t_1}\, $ and there exists another closed disk attached to this point whose interior has empty intersection with $\, sp ( x )\, $. Then it is easy to see that $\, z_1\in\delta sp ( x ) = {\delta }_e sp ( x )\, $ and there exists a holomorphic function in some neighbourhood of 
$\, sp ( x )\, $ taking its maximal absolute value in $\, z_1\, $ and strictly smaller in any other point of 
$\, sp ( x )\, $. But since the restriction map $\, {\mathfrak H}_e ( sp ( x ) ) = \mathfrak H ( sp ( x ) ) \rightarrow \mathfrak H ( X )\, $ is an isometric inclusion this leads to a contradiction. By a similar argument one gets that $\, sp ( x )\, $ cannot have any holes, i.e. its intersection with the holes of 
$\, X\, $ must be dense so by compactness any point confined by the outer boundary of $\, X\, $ must be contained in $\, sp ( x )\, $ showing that $\, sp ( x )\, $ is precisely the disjoint union of simply connected domains obtained by filling up the holes in $\, X\, $.  By approximating an arbitrary disjoint union of simply connected compact sets $\, Y\subset\mathbb C\, $ by domains of this type one easily deduces the result $\, \mathfrak H ( Y ) = {\mathfrak H}_e ( Y )\, $ for such a subset reducing assertion III to assertion II for this type. We interrupt the continuing argument to record this important observation.
\par\bigskip\noindent
{\bf Corollary.}\quad Let $\, Y\subseteq \mathbb C\, $ be a finite disjoint union of simply connected compact sets $\, \{ Y_j \}\, $. Then any holomorphic function defined in some open neighbourhood of 
$\, Y\, $ can be uniformly approximated on $\, Y\, $ by holomorphic polynomials\qed
\par\bigskip\noindent
Assertion II will follow in general from the inductive limit argument below which then shows that 
$\, {\widehat{\mathfrak H}}_e ( X ) \simeq {\widehat{\mathfrak H}}_e ( sp ( x ) ) \simeq C ( \delta sp ( x ) )
= C ( {\delta }_e X )\, $ proving assertion III for $\, X\, $.
The argument suffices to prove assertions I--III for any space $\, X\, $ which is the finite union of closed disks say of the same diameter (modulo the inductive limit argument below). Any compact subspace of the complex plane may be approximated by a sequence of these with common diameter of the disks shrinking to zero as $\, n\to\infty\, $.
\par\noindent
We now proceed with the inductive limit argument to obtain the general case. Let 
$\, \{ X_n \}\, ,\, X_nÊ\supseteq X_{n+1}\, $ for all $\, n\, $, be a decreasing sequence of compact subsets of the complex plane, each single space being the closure of a domain satisfying assertion I (this is sufficient for our purposes although the argument works the same assuming assertion II resp. assertion III for the $\, \{ X_n \}\, $), and let 
$\, X = {\cap }_{n=1}^\infty\, X_n\, $ be the limit space, so that by assumption
$$ \widehat{\mathfrak H} ( X )\> \simeq\> \lim_{\buildrel\rightarrow\over n}\, 
{\widehat{\mathfrak H}}_c ( X_n )\> \simeq\> \lim_{\buildrel\rightarrow\over n}\, C ( \delta X_n ) \> . $$ 
for the natural restriction maps of harmonic functions. For a sequence of operator spaces 
$\, \{ {\mathfrak X}_n \}\, $ let $\, ( \Pi / \oplus )_n {\mathfrak X}_n = {\Pi }_n {\mathfrak X}_n / \oplus_n {\mathfrak X}_n\, $ denote their outer direct product (on dividing their operator space direct product by the direct sum). As with the inverse limit one may also represent the direct limit of a sequence of spaces as a subspace of the outer direct product 
$$ \lim_{\buildrel\rightarrow\over n}\, {\widehat{\mathfrak H}}_c ( X_n )\> \subseteq\> 
( \Pi / \oplus )_n\, {\widehat{\mathfrak H}}_c ( X_n )\>\simeq\> ( \Pi / \oplus )_n\, C ( \delta X_n ) $$
by considering sequences that are compatible for the restriction maps. One notes that restricted to holomorphic functions the corresponding inclusion 
$$ \mathfrak H ( X )\> \subseteq\> ( \Pi / \oplus )_n {\mathfrak H}_c ( X_n ) $$
is even multiplicative for the structure of a unital operator algebra on the outer direct product equipped with the corresponding quotient matrix norm of the quotient space given by the direct product modulo the direct sum, both of which are immediately seen to be unital operator algebras. In particular the outer direct product of holomorphic function spaces is a unital operator space with enveloping operator system given by the outer direct product of the corresponding enveloping operator systems of the factors so the quotient norm is the same as for  the larger quotient (which by assumption on the $\, X_n\, $ defines a $C^*$-norm). Let 
$\, r_m : {\widehat{\mathfrak H}}_c ( X_m ) \rightarrow \widehat{\mathfrak H} ( X )\, $ be the restriction map. Choose an increasing sequence $\, \{ {\mathfrak F}_n {\}}_n\, $ of finitedimensional operator subsystems of $\, \widehat{\mathfrak H} ( X )\, $ with dense union.
One may assume by omitting certain of the 
approximating spaces and renumeration if necessary that for each $\, n\, $ there exists a finitedimensional operator subsystem $\, {\mathfrak F}^n\subseteq {\widehat{\mathfrak H}}_c ( X_n )\, $ 
such that the restriction map takes $\, {\mathfrak F}^n\, $ bijectively onto $\, {\mathfrak F}_n\, $ and is bounded below by $\, 1 - {\epsilon }_n\, $ on this subspace for a suitable sequence of positive numbers 
$\,\{ {\epsilon }_n \}\, $ with $\, {\epsilon }_n < 1 / (n + 1)\, $.
Then there exists from Lemma 1 of \cite{Ha2} by finite-dimensionality of the domain and commutativity of the range an extension of the map 
$\, {\mathfrak F}_n \rightarrow {\mathfrak F}^n\subseteq C ( \delta X_n )\, $ to a selfadjoint linear map 
$$ s_n :  C ( \delta X )\> \longrightarrow\> C ( \delta X_n ) $$ 
which is completely bounded by $\, ( n + 1 ) / n\, $. This leads for each $\, f\in C ( \delta X )\, $ to a sequence of elements 
$\, ( h_n ( f ) )_n = ( ( r_n\circ s_n ) ( f ) )_n \in \prod_n \widehat{\mathfrak H} ( X )\, $ by composition with the restriction maps $\, r_n\, $. We would like to have that this sequence converges uniformly to $\, f\, $ but  we can only prove pointwise convergence which however turns out sufficient for our purposes. Put 
$$ s = ( \Pi / \oplus )_n\, s_n\> :\> C ( \delta X )\> \longrightarrow\> ( \Pi / \oplus )_n\, C ( \delta X_n ) $$
and check that $\, s\, $ is a completely contractive extension of the natural inclusion of 
$\, \widehat{\mathfrak H} ( X )\, $ as above. Composing $\, s\, $ with the restriction map 
$\, r = (\Pi / \oplus )_n\, r_n\, $ gives a map 
$$ r\circ s :\> C ( \delta X )\> \longrightarrow ( \Pi / \oplus )_n \widehat{\mathfrak H} ( X )  $$
which embeds the subspace $\, \widehat{\mathfrak H} ( X )\, $ into the outer direct product by (images of) diagonal sequences. Let $\, \omega\, $ denote a free ultrafilter on $\, \mathbb N\, $ and consider the induced $*$-homomorphism 
$$ q_{\omega } :\> ( \Pi / \oplus )_n\,  l^\infty  ( \delta X ) \longrightarrow l^\infty  ( \delta X )  $$
by taking the pointwise limit along $\,\omega\, $.  Applying 
$\, q_{\omega }\, $ to an element of the form 
$\, ( r\circ s ) ( h )\, $ with $\, h\in \widehat{\mathfrak H} ( X )\, $ will give back the original element 
$\, h\, $ since the pointwise limit of a convergent sequence along $\,\omega\, $ is just the same as the ordinary limit. Let $\, t\in \delta X\, $ be any point and consider the linear functional 
$\, {\phi }_t\, $ on $\, C ( \delta X )\, $ given by the difference 
$\,  {\phi }_t ( f ) = ( ( q_{\omega }\circ r\circ s ) ( f ) )( t ) - f ( t )\, $. 
We claim that it is positive. Any strictly positive element $\, f\geq \epsilon {\bf 1}\, ,\, \epsilon > 0\, $ in  
$\, C ( \delta X )_+\, $ is the supremum of positive elements of the form $\, g = h\overline h\, $ with 
$\, h\in \mathfrak H ( X )\, $.  To see this it is sufficient to find for given such $\, f\, $ and $\, t\in\delta X \, $ a function $\, g\leq f\, $ of the prescribed form with $\, g ( t ) = f ( t )\, $. Without loss of generality one may assume that $\, f ( t ) = 1\, $. Choose some holomorphic function $\, k\, $ of norm one in $\, X\, $ taking its maximum in $\, t\, $ and such that $\, k ( s )\overline k ( s ) < 1\, $ is strictly smaller than $\, 1\, $ at any other point $\, s\neq t\, $ which exists since $\, \delta X\, $ is the Shilov boundary of 
$\,\mathfrak H ( X )\, $.  Taking a suitable power of $\, h = k^n\, $ one can achieve that 
$\, g = h\overline h \leq f\, $. Then 
$$ {\phi }_t ( f ) = ( e_t\circ q_{\omega }\circ r\circ s ) ( f ) - e_t ( f )\>\geq\> 
\sup\, \{ ( e_t\circ q_{\omega }\circ r\circ s ) ( h\overline h )\,\vert\, h\overline h \leq f\, ,\, h\in \mathfrak H ( X ) \} - e_t ( f ) $$
$$\quad \>\geq\> \sup\, \{ h ( t ) \overline h ( t )\,\vert\, h\overline h\leq f\, ,\, f\in\mathfrak H ( X ) \} - f ( t )
\> =\> 0 $$ 
since the image of a supremum under a positive map is larger or equal than the supremum of its images plus the Schwarz inequality for the positive map $\, e_t\circ q_{\omega }\circ r\circ s\, $ with $\, e_t\, $ evaluation at $\, t\, $. Therefore $\, {\phi }_t ( f )\, $ is positive for all $\, f \geq \epsilon {\bf 1}\, $ with arbitrary $\, \epsilon > 0\, $. By continuity $\, {\phi }_t\, $ is positive. Since $\, {\phi }_t ( 1 ) = 0\, $ one gets $\, {\phi }_t\equiv 0\, $. Therefore
$\, q_{\omega }\circ r\circ s \equiv id_{C ( \delta X )}\, $ for each and every free ultrafilter $\,\omega\, $ which can be the case only if the sequence $\, h_n ( f ) = r_n\circ s_n ( f )\, $ converges pointwise, i.e. weakly towards $\, f\, $. By a wellknown result of Banach space theory, cf. \cite{Pe} 2.4.8, there exists a sequence of convex combinations of the elements $\, \{ h_n ( f ) \}\, $ converging in norm towards $\, f\, $. This implies 
$\, f\in \widehat{\mathfrak H} ( X )\, $ and hence $\, C ( \delta X ) = \widehat{\mathfrak H} ( X )\, $. The proof of assertion III in the general case is just a copy of the argument above, replacing 
$\, C ( \delta X )\, $ with $\, C ( {\delta }_e X )\, $ etc., so assertions II and III are valid for an arbitrary compact subset $\, X\subset \mathbb C\, $. 
\par\medskip\noindent  
Given an element $\, x\, $ in a unital operator algebra $\, A\, $ let $\, sp ( x )\, $ denote its spectrum in the subalgebra $\, A_x\, $ generated by $\, x\, $ and the identity element, and 
$\, sp_A ( x )\, $ the spectrum of $\, x\, $ in $\, A\, $. Using holomorphic function calculus one may define the subalgebra $\,  A^{\circ }_x\subseteq A\, $ by considering the closure of the subalgebra generated by the union of all images of holomorphic functions $\, f\in \mathfrak H ( U )\, $ with 
$\, U\, $ an open neighbourhood of $\, sp_A ( x )\, $ in $\, A\, $ under the corresponding holomorphic function calculus. If $\, A\, $ embeds into a $C^*$-algebra such that $\, x\, $ is normal then one checks that each holomorphic function calculus 
$$ {\omega }_{U , x} : \mathfrak H ( U ) \longrightarrow A_x^{\circ } $$
is contractive, hence one obtains a completely isometric (bijective) holomorphic function calculus 
$$ {\omega }_{x , A} : \mathfrak H ( sp_A ( x ) )\> \longrightarrow\> A_x^{\circ }\> $$
in the limit. One may call $\, A_x^{\circ }\, $ {\it the subalgebra generated by the spectrum of $\, x\, $ in 
$\, A\, $}. An element $\, x\, $ will be called {\it weakly ($C^*$)-normal} iff for each open neighbourhood 
$\, U\, $ of $\, sp_A ( x )\, $ the holomorphic function calculus 
$\, {\omega }_{U , x}, $ of functions holomorphic in $\, U\, $ into $\, A\, $ is contractive and {\it locally ($C^*$)-normal} iff for each open neighbourhood $\, V\, $ of $\, sp ( x )\, $ the holomorphic function calculus $\, {\omega }_{V , x}\, $ is contractive. The prefix $C^*$ is used to distinguish this notion from the notion of normal element in a super $C^*$-algebra with respect to the superinvolution and can be dropped if no confusion can arise. These definitions also make sense in case of a Banach algebra. The following result however is meaningful only for operator algebras.
\par\bigskip\noindent
{\bf Theorem}\quad (Harmonic function calculus) 
If $\, X\subseteq \mathbb C\, $ is any compact subset and $\, \delta X\, $, resp. $\, {\delta }_e X\, $, denotes the Shilov boundary of 
$\, \mathfrak H ( X )\, $, resp. $\, {\mathfrak H}_e ( X )\, $, then 
$\, C ( \delta X )\, $ is equal to the enveloping operator system of harmonic functions $\, \widehat{\mathfrak H} ( X )\, $ (in its injective envelope) and 
$\, C ( {\delta }_e X )\, $ is equal to the enveloping operator system $\, {\widehat{\mathfrak H}}_e ( X )\, $. Given a weakly normal element $\, x\in A\, $ of a unital operator algebra (note that this implies that the spectral radius of $\, x\, $ is equal to its norm)
the holomorphic function calculus 
$$    {\omega }_{x , A} :  \mathfrak H ( sp_A ( x )  )\> \longrightarrow\> A $$
sending the function $\, id_{sp_A ( x )} \in \mathfrak H ( sp_A ( x ) )\, $ to $\, x\, $ and $\, 1\, $ to $\, 1\, $ is well defined and completely isometric, hence extends uniquely to a completely (positive) isometric map on the enveloping operator systems  
$$ {\widehat \omega}_{x , A} :\> C ( \delta sp_A ( x ) )\simeq 
\widehat{\mathfrak H} ( sp_A ( x ) )\> \longrightarrow\> {\widehat A}_x^{\circ }\>\subseteq \widehat A \> . 
$$
In the same manner given a locally normal element $\, x\in A\, $ the holomorphic function calculus 
$$ {\omega }_x : \mathfrak H ( sp ( x ) ) = {\mathfrak H}_e ( sp_A ( x ) )\> \longrightarrow\>  A $$
is completely isometric, hence extends uniquely to a completely (positive) isometric map on the enveloping operator systems
$$ {\widehat\omega }_x :\>  C ( \delta sp ( x ) ) = C ( {\delta }_e sp_A ( x ) )\> \longrightarrow\> 
{\widehat A}_x\>\subseteq\> \widehat A\> . $$
For an arbitrary element $\, x\in A\, $ the holomorphic function calculus in a given neighbourhood of 
$\, sp_A ( x )\, $ is completely contractive whenever it is contractive (for example if choosing the neighbourhood so large as to contain the disk of radius 
$\,\Vert x \Vert\, $ centered in the origin).
\par\bigskip\noindent
{\it Proof.}\quad The first part has been proved above. Assume then that $\, x\, $ is weakly normal. 
This implies that the holomorphic function calculus $\, {\omega }_{x , A}\, $ is well defined unital and contractive so it extends to a positive map $\, {\widehat\omega }_{x , A}\, $ on the enveloping operator systems by Proposition 2.12. of \cite{Pa}. Since the domain of the extension is a commutative $C^*$-algebra the extension is completely positive from Theorem 5.1.5 of \cite{E-R}. Then the holomorphic function calculus must be completely contractive as well. On the other hand composition with the completely contractive Gelfand transform 
$$ {\mathfrak G}_{x , A} :\> A_x^{\circ }\> \longrightarrow \mathfrak H ( sp_A ( x ) ) $$
gives an inverse for $\, {\omega }_{x , A}\, $ so that both maps must be completely isometric. The proof in case of a locally normal element is just the same\qed
 \par\bigskip\noindent
By aid of this theorem we may define a (weak) superpositive absolute value $\, \vert x {\vert }_{s , w}\, $ for any weakly $C^*$-normal hermitian element $\, x\, $ in a super-$C^*$-algebra which is an element of its enveloping graded operator system, by considering the image of the graded absolute value of 
$\, id_z\, $  in $\, C ( \delta sp_A ( x ) )\simeq \widehat{\mathfrak H} ( sp_A ( x ) )\, $ under the corresponding harmonic function calculus. For a locally $C^*$-normal hermitian element one defines the local superpositive absolute value $\, \vert x {\vert }_{s , l}\, $ by considering the image of the graded absolute value of $\, id_z\, $ in $\, C ( \delta sp ( x ) )\simeq \widehat{\mathfrak H} ( sp ( x ) )\, $ under harmonic function calculus. 
The suffix $\, s\, $ is chosen to distinguish the superpositive absolute value from the ordinary absolute value (which can also be defined in the enveloping graded operator system by harmonic function calculus) and is termed $\, \vert x {\vert }_w\, $ for a weakly $C^*$-normal element resp. 
$\, \vert x {\vert }_l\, $ in case of a locally $C^*$-normal element.
In $\, C ( \delta sp ( x ) )\, $ the element 
$\, \vert x {\vert }_{s , l}\, $ is given by the same formula written out above in case of a closed disk (by restriction to $\,\delta sp ( x ) \, $). In particular $\, \Vert \vert x {\vert }_{s , l} \Vert = \Vert x \Vert\, $.
\par\bigskip\noindent
The Dirichlet problem applies to more general domains in the Euclidean space $\, {\mathbb R}^n\, $ where the connection of harmonic functions as enveloping operator system of an algebra of holomorphic functions is lost. One may define however an operator algebra consisting of complexvalued harmonic functions in any dimension, at least for prime dimensions. We give the definition for 
$\, n = 3\, $, the generalization to higher dimensions is immediate. Namely, consider the set of complexvalued harmonic functions defined in some neighbourhood of a fixed compact subspace $\, X \subset {\mathbb R}^3\, $ such that each element takes the form 
$$ \mathfrak x ( x , y , z ) \> =\> f ( x , y , z )\, +\, \sqrt\omega\, g ( x , y , z )\, +\, \sqrt{{\omega }^2}\, h ( x , y , z ) $$
where $\, f\, ,\, g\, ,\, h\, $ are realvalued twice continuously differentiable functions and $\,\omega\in\mathbb C\, $ is a third root of unity, subject to the partial differential equations
$$  \partial f / \partial x\> =\> - \partial g / \partial y\> =\> \partial h / \partial z\> ,\quad \partial f / \partial y\> =\> - \partial g / \partial z\> =\> \partial h / \partial x\> , $$
$$ \partial f / \partial z\> =\> - \partial g / \partial x\> =\> \partial h / \partial y\> . $$ 
One checks that $\, \mathfrak x\, $ is harmonic and the product of two such elements can be written in the same form so that this set of functions defines an algebra whose uniform closure in $\, C ( X )\, $ is a uniform operator algebra which may be called the algebra of {\it triholomorphic functions} in $\, X\, $. It should be interesting to determine its Shilov boundary and examine if a similar correspondence to the class of all harmonic functions over $\, X\, $ can be established as with holomorphic functions.
\par

\end{document}